\documentclass[12pt]{article}

\input{amssym.def}
\input{amssym}
\usepackage{eufrak}

\input{amssym.def}
\input{amssym}

\def\theequation{\thesection.\arabic{equation}}
\makeatletter \@addtoreset{equation}{section} \makeatother
\def\nn{\nonumber}
\def\be{\begin{equation}}
\def\ee{\end{equation}}
\def\ba{\begin{eqnarray}}
\def\ea{\end{eqnarray}}

\def\tr{{\rm tr}}
\def\Tr{{\rm Tr}}

\def\tbl#1#2{{\ifmmode
\Bigl\{\!\!\begin{array}{c}
\scriptstyle #1\\[-2pt]
\raisebox{2pt}{$\scriptstyle #2$} \end{array}\!\!\Bigr\} \else
\raisebox{2pt}{\scriptsize$\left\{\!\!\!\begin{array}{c}
#1\\[-1pt] #2 \end{array}\!\!\!\right\}$}\fi}}

\newlength{\coefkappa}
\begin{document}

\def\h{{\hbar}}
\def\Cas{{\rm Cas}}

\def\k{\mathbf{k}}
\def\qq{q^{-1}}
\def\qqq{q^{-2}}
\def\gq{\gggg_{q}}

\def\uq{U_q(sl(2))}
\def\Uq{U_q(sl(n))}
\def\Uqq{U_q(\gggg)}

\def\C{{\Bbb C}}
\def\R{{\Bbb R}}
\def\K{{\Bbb K}}

\def\De{\Delta}
\def\de{\delta}

\def\ve{\varepsilon}
\def\ep{\epsilon}
\def\al{{\alpha}}
\def\aal{{\alpha}^{-1}}
\def\ot{\otimes}
\def\om{\omega}
\def\Om{{\Omega}}
\def\End{{\rm End\, }}
\def\Tr{{\rm Tr}}
\def\Sym{{\rm Sym\, }}
\def\Gr{{\rm Gr}}
\def\oi{\overline{i}}
\def\oj{\overline{j}}

\def\dim{{\rm dim}}
\def\sdim{{\rm sdim}}
\def\vv{V^{\ot 2}}
\def\vl{V_{\la}}
\def\la{{\lambda}}
\def\span{\rm span}
\def\lhq{{\cal L}(\h,q)}
\def\lqh{{\cal L}(\h,q)}
\def\slqh{{\cal SL}(\h,q)}
\def\slqh{{\cal SL}(\h,q)}
\def\lqo{{\cal L}(1,q)}
\def\lq{{\cal L}(q)}
\def\L{{\cal L}}
\def\SL{{\cal SL}}
\def\lh{{\cal L}_\h}
\def\F{{\cal F}}
\def\A{{\cal A}}
\def\aq{{\cal A}(q)}
\def\aqh{{\cal A}(\h,q)}
\def\ahq{{\cal A}(\h,q)}
\def\smn{{\sigma^{m,n}}}
\def\RR{{{\textsf{R}}}}
\def\Omn{\K_q[{\cal O}_{\mu, \nu}]}
\def\O{{\cal O}}
\def\M{{\cal M}}
\def\E{{\cal E}}
\def\gg{{\frak g}}
\def\gh{{\frak g}_\h}

\def\bea{\begin{eqnarray}}
\def\eea{\end{eqnarray}}
\def\be{\begin{equation}}
\def\ee{\end{equation}}
\def\nn{\nonumber}

\makeatletter
\renewcommand{\theequation}{{\thesection}.{\arabic{equation}}}
\@addtoreset{equation}{section} \makeatother

\title{Quantization of pencils with a $gl$-type Poisson center and braided geometry}
\author{
\rule{0pt}{7mm} Dimitri
Gurevich\thanks{gurevich@univ-valenciennes.fr}\\
{\small\it LAMAV, Universit\'e de Valenciennes,
59313 Valenciennes, France}\\
\rule{0pt}{7mm} Pavel Saponov\thanks{Pavel.Saponov@ihep.ru}\\
{\small\it Division of Theoretical Physics, IHEP, 142284 Protvino,
Russia} }

\maketitle

\begin{abstract}
In the algebra $\Sym(gl(m))$ we consider  Poisson pencils generated by the linear Poisson-Lie bracket $\{\,,\,\}_{gl(m)}$
and that corresponding to the so-called Reflection Equation Algebra. Each bracket of such a pencil has the Poisson center
coinciding with that of the bracket $\{\, ,\,\}_{gl(m)}$. Consequently, any bracket from this pencil can be restricted  to a
generic $GL(m)$-orbit $\O\subset gl(m)^*$. Quantization of such a restricted bracket can be done in the frameworks of
braided affine geometry. In the paper we consider these Poisson structures, their super-analogs as well as their quantum
(braided) counterparts. Also, we exhibit some detailed examples.
\end{abstract}

{\bf AMS Mathematics Subject Classification, 2010:} 81R60, 81S99


\section{Introduction}
\label{sec:int}

In the paper we deal with certain Poisson pencils defined in the algebras $\K[gl(m)^*]\cong \Sym(gl(m))$ and their
super-analogs. Hereafter, $\K$ is the ground field, namely $\C$ or $\R$. The Poisson center\footnote{By the Poisson center
(called  in the sequel the {\em center}) we mean the set of functions $f$ Poisson commuting with any other function $g$:
$\{f,g\}=0$.} of each bracket from such a pencil coincides with that of the linear Poisson-Lie bracket $\{\,,\,\}_{gl(m)}$ (or
that of $\{\,,\,\}_{gl(m|n)}$ in a super-case) coming in the pencil. It is in this sense that we speak about pencils with the
$gl$-type center. Consequently, any  bracket from such a pencil can be restricted to an arbitrary generic $GL(m)$-orbit
$\O\subset gl(m)^*$ (or its super-analog).

The simplest example of pencils possessing this property is the following one (below it will be treated as a restriction of a
pencil defined in the algebra $\Sym(gl(2))$). Let
$$
\{\,,\,\}_{sl(2)}:\;\Sym(sl(2))^{\ot 2}\to \Sym(sl(2))
$$
be the linear Poisson-Lie bracket corresponding to the $sl(2)$ Lie structure and defined in the standard basis $\{x,h,y\}$
of the algebra $\Sym(sl(2))$ as follows
\be
\{x,y\}_{sl(2)}=h, \quad \{h,x\}_{sl(2)}=2x, \quad \{h,y\}_{sl(2)}=-2y.
\label{zero}
\ee
The Poisson center of the bracket (\ref{zero}) consists of functions $p(xy+yx+{{h^2}\over{2}})$ where $p$ is a polynomial
in one variable.

Also, consider the following quadratic Poisson bracket
\be
\{x,y\}'=h^2, \quad \{h,x\}'=2xh, \quad \{h,y\}'=-2yh.
\label{odin}
\ee
It is easy to see that these two brackets are compatible and any bracket from the corresponding pencil
\be
\{\,,\,\}_{a,b}=a\,\{\,,\,\}_{sl(2)}+b\, \{\,,\,\}'
\label{penc}
\ee
possesses the same center as the bracket $\{\, ,\,\}_{sl(2)}$ does. Consequently, this Poisson pencil can be restricted to any
variety defined by the equation $xy+yx+{{h^2}\over{2}}=C\not=0$.

In virtue of the famous  Kontsevich result \cite{K} any Poisson structure on a smooth variety $\M$ can be quantized by
deformation quantization means. Namely, there exists a new associative product in the commutative algebra
$\K[\M][[\h]]$ (where $\h$ is a quantization parameter), satisfying the so-called correspondence principle.
Consequently, each individual bracket from the pencil (\ref{penc}) or its restriction can be quantized in this sense.
However, in general it is not clear what are relations between quantum algebras arising from the pencil (\ref{penc}) and
those arising from its restrictions. Namely, whether the latter quantum algebras can be realized as some quotients of the
former ones. In order to answer this question we have to describe the quantum analog of the center of the pencil (\ref{penc}).

Fortunately, the Poisson pencil (\ref{zero}) can be explicitly quantized. As a result, we get a family of quantum algebras
depending on  two parameters (which can be specialized to  numbers since our quantization is not formal). Moreover, the
center of any such an algebra can be easily described. Namely, each center is also generated by a quadratic element but it
is not symmetric any more with respect to permutation of the factors in its summands  and it cannot be written in a symmetric
form. Consequently, the pairing defined on the space $\L={\span}_\K(x,h,y)$ via the matrix inverse to that formed by the
coefficients of this quadratic element is not symmetric either. This is a hint that the resulting quantum algebra is related to a
braiding, i.e. to a solution of the quantum Yang-Baxter equation (see section \ref{sec:3}). Though in section \ref{sec:2} we
quantize the pencil (\ref{zero}) by a direct and somewhat elementary method, it can be also done by a general method based
on the so-called $R$-matrix technique. This technique enables us to introduce a quantum (braided) trace (or $R$-trace) which
plays a crucial role in describing the center of the algebras (super-algebras included) arising from pencils similar to (\ref{penc}).

The $R$-trace is an ingredient of the braided geometry considered in \cite{GS1}. Other ingredients are braided Lie
algebras, braided vector fields (which are not considered in this paper), and braided affine varieties.
Braided varieties, we are dealing with, are deformations of generic $GL(m)$-orbits in $gl(m)^*$ (or their super-analogs).
They are in a sense regular varieties since for such a variety there exists a projective module playing the role of the cotangent
vector bundle in the frameworks of the Serre approach (see \cite{GS2}). Thus, by quantizing $gl(m)$ generalization of the
Poisson pencil (\ref{penc}) we get a braided (i.e. related to a braiding) deformation of the enveloping algebra $U(gl(m))$. Whereas the quantum counterpart of such a pencil restricted to a generic orbit in $gl(m)^*$ can be realized as a quotient
of this "braided enveloping algebra". Finally, this quotient is treated to be a braided generic orbit.

Note that Poisson pencils analogous to (\ref{penc}) exist on super-algebras $\Sym(gl(m|n))$ too.
Their quantization gives rise to braided algebras with similar properties. The main goal of this paper is to describe these
Poisson structures, their restrictions to generic super-orbits in $gl(m|n)^*$ and their quantum counterparts. A crucial
role in our construction is played by a quantum (braided) version of the Cayley-Hamilton (CH) identity valid for a generating
matrix $L$ of the Reflection Equation Algebra (REA) found in \cite{GPS1}. This identity enables us to define eigenvalues of
the matrix $L$. In terms of these eigenvalues we introduce a criterium of regularity of super-orbits and their braided counterparts.

The paper is organized as follows. In the next section we quantize the Poisson pencil (\ref{penc}) without using any technique
related to braidings. Nevertheless, we arrive to an algebra  which is in a sense a braided algebra. In section \ref{sec:3} we
consider a family of such braided algebras deforming  super-commutative algebras $\Sym(gl(m|n))$ and exhibit the
corresponding super-Poisson structures. In section \ref{sec:4} we consider restrictions of such Poisson structures to generic
super-orbits in $gl(m|n)^*$. Their quantum counterparts are braided generic orbits. In section \ref{sec:5} we consider two
low-dimensional examples (the first of them is just our basic example but treated in the frameworks of  braided geometry).
In the last section we list a few open problems related to a more general class of Poisson pencils with $gl$-type center.
\medskip

\noindent
{\bf Acknowledgement.} The work of one of the authors (P.S.) was partially supported by
the RFBR grant 08-01-00392-a and the joint RFBR and DFG grant 08-01-91953. The work of D.G. and
P.S. was partially supported by the joint RFBR and CNRS grant 09-01-93107.

\section{Basic example}
\label{sec:2}

Let us consider the  algebra $\K[sl(2)^*]\cong \Sym(sl(2))$  endowed  with the Poisson pencil (\ref{penc}).
Observe that the Poisson bracket (\ref{odin}) coming in this pencil is quadratic and  differs by the factor $h$ from the linear bracket (\ref{zero}).
Furthermore, it can be easily seen that
 the function $\Cas=xy+yx+{{h^2}\over{2}}$ is central for any Poisson bracket from this pencil:
$$
\{\Cas, f\}_{a,b}=0\quad \forall f\in \Sym(sl(2)).
$$

We treat the enveloping algebra $U(sl(2)_\h)$ of the Lie algebra\footnote{The notation $\gh$ means that we
introduce the factor $\h\in \K$ in the  Lie bracket of the Lie algebra $\gg$.} $sl(2)_\h$ to be a quantum counterpart
of the Poisson algebra $\Sym(sl(2))$ with the Poisson-Lie bracket $\{\,,\,\}_{sl(2)}$. Our immediate aim is to quantize any
bracket from the pencil  (\ref{penc}) (in fact, we simultaneously quantize  the whole Poisson pencil).

First, we quantize the bracket $\{\,,\,\}'$ alone. Consider an associative algebra generated by
three elements $x,h,y$ subject to the relations
\be
hx-xh=\nu\, (a\, hx+b\,xh),\quad hy-yh=-\nu\, (c\,hy+d\,yh),\quad xy-yx=\nu\, h^2,
\label{quant}
\ee
where $a$, $b$, $c$ and $d$ are parameters subject to the constraint $a+b=c+d=2$. The quantization parameter
$\nu$ (as well as all  parameters below) can be  specialized to a number from the ground field $\K$.

The main feature of a quantization of the algebras in question is that it should give rise to quantum objects with a good
deformation property in the following sense. Let a quadratic algebra $\A(\nu)$ depends on a parameter $\nu$ and at the "classical limit" $\nu\to 0$
it turns into the symmetric algebra of a space $V$ : $\A(0) =\Sym(V)$. We say that $\A(\nu)$ possesses a {\em good
deformation property} if $\dim\, \A(\nu)^k=\dim\, \Sym^k(V) $  for any $k\ge 0$ and a generic $\nu$. Here the superscript
$k$ stands for the $k$-th degree homogeneous component. If moreover, $\A(\nu, \h)$ is a quadratic-linear algebra
such that $\A(\nu, 0)=\A(\nu)$, we say that it has a {\em good deformation property} if $\Gr \A(\nu, \h)=\A(\nu)$
where $\Gr$ stands for the associated graded algebra. In the same sense we speak about the good deformation property
of algebras close to super-algebras $\Sym(gl(m|n))$.

In what follows we additionally assume $a=d$, $b=c$ in relations (\ref{quant}). Otherwise, as can be shown (see footnote
\ref{ftn:3}), the corresponding quotient algebra is not a quantum object,  i.e. it does not possess the good deformation property.
Under this condition we can rewrite relations (\ref{quant}) in the form
\be
q^2\, hx-xh=0,\quad q^2 yh-hy=0,\quad xy-yx-\nu\, h^2=0,\,\,{\rm where}\,\,  q^2 = {{1-a\nu}\over{1+b\nu}}.
\label{quant1}
\ee

In order to show that the algebra defined by relations (\ref{quant1}) is indeed a quantization of the Poisson algebra
$Sym(sl(2))$ with the Poisson bracket $\{\,,\,\}'$, we fix the family of elements $\{e_{k,l,m}=x^k y^l h^m, \;k,l, m=0,1,2,...\}$
in it. Then we have to show that this family is a basis of the algebra in question (an analog of the Poncar\'e-Birhoff-Witt theorem). To this end we have to check that the products $xe_{k,l,m}$, $ye_{k,l,m}$, and $he_{k,l,m}$ can be expressed as
a linear combinations of the elements $\{e_{k,l,m}\}$. Besides, we should verify, that relations
\be
(q^2hx- xh)\,  e_{k,l,m}=0,\quad (q^2 yh-hy)\, e_{k,l,m}=0,\quad (xy-yx-\nu\, h^2) \, e_{k,l,m}=0
\label{prop}
\ee
do not lead to any dependencies among the elements $e_{k,l,m}$ for all $k,l,m$. Details are left to the reader.

Note that another way of verifying the good deformation property of such type algebras is based on some special projectors
\cite{GPS2} (see also the next section).

Thus, we have got a family of quantum algebras depending on a value of $a$. Nevertheless, for $\nu\not=0$ all these
quantizations are equivalent (over $\C$). It can be shown by rescaling the generator $h$. So, we set $a=b=1$. Then  the
relations between generators become
\be
q^2hx-xh=0,\quad q^2 y h-hy=0,\quad (q^2+1)(x y-y x)+ (q^2-1) h^2=0
\label{sys}
\ee
(i.e. $q^2={{1-\nu}\over{1+\nu}}$ or equivalently, $\nu={{1-q^2}\over{1+q^2}}$).

Denote $\aq$ the algebra generated by the space  $\A={\span}_\K(x, h, y)$ where the generators are subject to the relations
(\ref{sys}). The algebra $\aq$ is in a sense  "$q$-symmetric" algebra of the space $\A$. Below, we explain the exact meaning
of this claim. Now, we pass to a quantization of the whole Poisson pencil (\ref{penc}).

To this end we look for numerical factors $A,B,C$ such that the algebra defined by relations
\be
q^2hx-xh=A\,x,\,\, q^2 y h-hy=B\,y,\,\,(q^2+1)(x y-y x)+ (q^2-1) h^2=C\,h
\label{quant2}
\ee
would have the good deformation property.

In order to find such factors we use the Jacobi identity in the form of  \cite{PP}. Let $I\subset \A^{\ot 2}$ be a subspace
spanned by the left hand side of (\ref{quant2}). Then the space $I\ot \A\bigcap \A\ot I$ is one-dimensional\footnote{\label{ftn:3}
If in formula (\ref{quant}) $a\not= d$ the space $I\ot \A \bigcap \A\ot I$ is trivial. This entails that $\dim(I\ot \A + \A\ot I)$ differs
from the classical one and therefore the algebra defined by (\ref{quant}) does not have the good deformation property.} and it
is spanned by the following element
$$
y(q^2 hx-xh)+ \qqq x(q^2 yh-hy)+h(xy-yx-\nu\, h^2)=
$$
$$
\qqq (q^2 hx-xh)y+(q^2 yh-hy)x+(xy-yx-\nu\, h^2)h.
$$
Using relations (\ref{quant2}) we reduce this equality to the form
$$
Ayx+\qqq Bxy+{{Ch^2}\over{1+q^2}}-(\qqq Axy+Byx +{{Ch^2}\over{1+q^2}})=0.
$$
According to the Jacobi condition from \cite{PP} the left hand side of this relation  must belong to $I$. Gathering similar
terms and applying the relations (\ref{quant2}), we come to the equality
$$
(A-B)\Big((1-q^{-2})xy +\frac{q^2-1}{q^2+1}\,h^2 -\frac{C}{q^2+1}\,h\Big) = 0.
$$
In case $q^2\not=1$ the only possible choice is $A=B$. For the factor $C$ there is no restriction.

We assume $C\not=0$, then by rescaling $x$ (or $y$) we can get\footnote{The
factor $2_q$ will be motivated in section \ref{sec:3}, see (\ref{motiv}).}
$$
C=A=B=2_q\h,
$$
where $\h$ is a new quantization parameter and the $q$-numbers are defined in the usual way
$$
k_q=\frac{q^k-q^{-k}}{q-q^{-1}},\quad k\in {\Bbb Z}.
$$

We denote $\aqh$ the algebra defined by the relations (\ref{quant2}) with $A=B=C=2_q \h$. So, if $\h=0$ this algebra turns into
$\aq$. Choosing the same basis $\{e_{k,l,m}\}$ in the algebra $\aqh$, $\h\not=0$, we can show that the property similar to
(\ref{prop}) (with linear  terms added) is still valid  and we conclude that the algebra $\aqh$ have the good deformation property
and it is a two parameter deformation of that $\Sym(sl(2))$.

In order to get  a quantization of one bracket from the pencil $\{\,,\,\}_{a,b}$ it suffices to bound the parameters of quantization
$q=exp(\al \mu)$, $\h=\beta \mu$, and to find the Poisson bracket corresponding to the parameter $\mu$.

Now, motivated by the fact that all brackets $\{\,,\,\}_{a,b}$ have the center generated by the element $\Cas = xy+yx+h^2/2$ we
want to find the center of the quantum algebra $\aqh$. It is not difficult to check that the element
$$
\Cas_q=\qq xy+ q yx +{{h^2}\over{2_q}}
$$
is central in the algebra $\aqh$ and therefore, so are all elements $p(\Cas_q)$ where $p$ is a polynomial in one variable.

The matter is that the element $\Cas_q$ is not symmetric and it cannot be written in a symmetric form. The pairing
$\A^{\ot 2}\to\K$ defined by the matrix inverse to a matrix composed from the coefficients of the element  $\Cas_q$
becomes\footnote{Note that this way of defining the pairing is motivated by identification of the spaces $V$ and $V^*$ in the
monoidal quasitensor rigid category generated by the space $V$ as described in \cite{GLS}.}
\be
\langle x,y \rangle=\qq,\quad \langle y,x \rangle=q,\quad \langle h, h\rangle=2_q.
\label{pair}
\ee
The pairing is not symmetric either. This is a hint that this quantum algebra can be related  to a braiding different from the usual
permutation operator. In section \ref{sec:5} we exhibit this relation after having considered the general case in sections \ref{sec:3}.
In section \ref{sec:4} we also consider general analogs of the {\em quantum (braided) hyperboloids} defined by the equation
$\Cas_q=const\not=0$.

\section{REA and corresponding Poisson pencils on super-spaces}
\label{sec:3}

In the previous section we considered an example of a Poisson pencil such that its brackets possess just the same center as the bracket
$\{\,,\,\}_{sl(2)}$ has.  Also, we quantized this pencil
without using any braiding. In this section we consider a general case which includes the previous example.
Our consideration also covers Poisson pencils on super-spaces $gl(m|n)^*$. Their quantization gives rise to algebras related to certain
braidings as well. However, our presentation goes in the opposite direction: we begin with the quantum objects called
(modified) Reflection Equation Algebras. Afterwards, we consider their Poisson counterparts. Each of these counterparts is a pencils comprising
the linear bracket $\{\,,\,\}_{gl(m|n)}$ and having the center of $gl$ type.
This property enables us to restrict the Poisson pencils in question to generic orbits in $gl(m|n)^*$.

Let us  consider a super-space $V=V_0\oplus V_1$ with $\dim V_0=m$ and $\dim V_1=n$. We call the ordered pair $(m|n)$
the {\em super-dimension} of the super-space $V$. Let $R\in\End(\vv)$ be a {\em Hecke symmetry}\footnote{Recall that
by a Hecke symmetry we mean a braiding $R\in\End(\vv)$, i.e. a solution of the quantum Yang-Baxter equation
$$
(R\ot I)(I\ot R)(R\ot I)=(I\ot R)(R\ot I)(I\ot R),
$$
subject to the second degree equation
$$
(qI-R)(\qq I+R)=0.
$$
We assume $q\in \K$ to be generic. In particular, this means that $q\not=0$ and $q^n\not=1$ for $n=2,3,4,...$}
defined as follows
\be
R=\sum_{1\leq i\leq m+n} (-1)^{\oi} \, q^{1-2\oi} e_{i}^i\ot e_{j}^j+\sum_{i\not=j} (-1)^{\oi\,\oj}
e_{i}^j\ot e_{j}^i +(q-\qq)\sum_{j>i} e_{i}^i\ot e_{j}^j,
\label{R-matr}
\ee
where $e_{i}^j$ stands for the $(m+n)\times (m+n)$ matrix with 1 at the intersection of the i-th row and j-th column and 0 otherwise and $\oi$ is the parity of $i$, i.e.
$$
i=0\quad {\rm  if} \quad 1\leq i\leq m \quad {\rm and} \quad i=1 \quad {\rm if}\quad m+1 \leq i \leq m+n.
$$
Note that for $q\to 1$ this braiding turns into a super-flip denoted in the sequel $\smn$.

Consider a unital associative algebra $\lqh$ generated by indeterminates $l_i^{j},\,\, 1\leq i,j\leq m+n$ subject to the
following multiplication rules
\be
R\,L_1\,R\,L_1-L_1\,R\,L_1\, R=\h(R\,L_1-L_1\,R),\quad L_1=L\ot I,\; L=\|l_i^{j}\|.
\label{RE}
\ee
 We call $L$ the {\em generating matrix} of the algebra $\lqh$. Here $\h\in \K$ and $q\in \K$ coming in the braiding $R$ are two
deformation parameters. The one parameter algebra $\lq:=\L(0,q)$ is called {\em Reflection Equation Algebra} (REA)
while the algebra $\lqh$ will be referred to as the {\em modified REA}.

Going back to the Hecke symmetry (\ref{R-matr}) note that it is {\em skew-invertible}. By definition, this means that there
exists an operator $\Psi \in \End(\vv)$ such that
$$
\Tr_2 R_{12} \Psi_{23}= \sigma_{13}.
$$
Here $\Tr_2$ stands for the (usual) trace applied to the operator product $R_{12} \Psi_{23}\in \End(V^{\ot 3})$  in the second
space and $\sigma_{13}$ is the usual flip transposing the first and third spaces in the tensor product $V^{\otimes 3} :=
V_1\otimes V_2\otimes V_3$.

Consider two operators $B:V\to V$ and $C:V\to V$ defined as follows
$$
B=\Tr_1 \Psi,\quad  C=\Tr_2 \Psi.
$$
Note the operators $B$ and $C$ are bound by the relation
\be
BC=q^{2(n-m)} I
\label{BC}
\ee
provided the Hecke symmetry $R$ is a deformation of the super-flip $\smn$ (and even in a more general setting discussed
in \cite{GPS2}). Consequently, the operators $B$ and $C$ are invertible.

These operators play a crucial role in defining an intrinsic trace $\Tr_R$ related to the braiding $R$. Namely, we put by
definition
$$
\Tr_R L^k:=\Tr (L^k C).
$$
We call the operation $\Tr_R$ the {\em $R$-trace}. The crucial property of the elements $\Tr_R L^k$ is that they
are central in the algebra $\lqh$. They are called {\em braided Casimir} elements. We are especially
interested in the  braided quadratic Casimir element $\Tr_R L^2$.

As  for the operator $B$, we use it for constructing a representation of the algebra $\lqh$. Namely, in the basis $\{x_i\}$ of
the space $V$ coordinated with the matrix form (\ref{R-matr}) of the above Hecke symmetry we set
$$
\pi(l_i^j)(x_k)=B_k^j x_i.
$$
Then the map
$$
\pi:\lqo\to \End(V)
$$
defines a representation of the algebra $\lqo$ (see \cite{GPS2}). Moreover, we get  an embedding $\L\to\End(V)$  where
$\L={\span}_{\K}(l_i^j)$ and consequently the family $\{l_i^j\}$ constitutes a basis of the space $\End(V)$.

We also need a numerical $R$-trace operator
\be
\tr_R: \End(V)\to \K,
\label{numtr}
\ee
which is a braided analog of the usual numerical trace. In the basis $\{l_i^j\}$ it has the form $\tr_R(l_i^j)=\de_i^j$.
Note that the above embedding $\L\to \End(V)$ also enables us to present the usual product
$$
\circ:\End(V)^{\ot 2}\to \End(V)
$$
in the basis $\{l_i^j\}$. Namely, we have $l_i^j\circ l_k^m=B_k^j\,l_i^m$. Consequently, we can define a pairing on the space
$\L$ by setting
\be
\langle\,,\,\rangle: \L^{\ot 2}\to \K,\quad \langle l_i^j, l_k^m\rangle= \tr_R (l_i^j\circ l_k^m)=B_k^j\,\de_i^m.
\label{pai}
\ee

Note that this pairing is non-degenerate on the space $\L$. Also, as follows from the relation (\ref{BC}), the matrix of this
paring is inverse (up to a factor) to the matrix of coefficients in the braided quadratic Casimir $\Tr_R L^2$.

In the case of the Hecke symmetry (\ref{R-matr}) the operator $C$  represented in the same basis $\{x_i\}$ of
the space $V$ has the form (see \cite{I})
$$
C_i^j=(-1)^{\oi} q^{2n+(-1)^{\oi}(2i-2m-1)} \de_i^j.
$$
The operator $B$ can be found from the relation (\ref{BC}).

We treat the  algebra $\lqh$ corresponding to the Hecke symmetry (\ref{R-matr}) to be a {\em braided} analog of the
enveloping algebra $U(gl(m|n)_\h)$. Now, we want to define an analog of the algebra $U(sl(m|n)_\h)$, provided that
$m\not=n$.

Let $\ell=\Tr_R L=\Tr(L\,C)$ be the linear braided Casimir element. Applying the numerical $R$-trace $\tr_R$ to this element
we have (see \cite{GPS2})
$$
\tr_R \ell=\Tr C=q^{m-n}(m-n)_q.
$$
This quantity vanishes iff $m=n$ (recall that $q$ is generic).

Assuming $m\not=n$, we put $\slqh=\lqh/\langle \ell \rangle$. The algebra $\slqh$ is generated by the space $\SL$ of
traceless elements (with respect to numerical $R$-trace $\tr_R$) of the space $\L$. Being restricted to the space $\SL$, the
pairing (\ref{pai}) is still non-degenerate. The algebra $\slqh$ is considered as a braided analog of the enveloping algebra $U(sl(m|n)_\h)$.

Now, consider  the Poisson structures corresponding to the  algebras $\lq$ and  $\lqh$, assuming the Hecke symmetry
entering their definition to be of the form (\ref{R-matr}). First, we consider the algebra $\lq$.
This algebra  has the good deformation property, i.e.
$\dim\, \lq^k=\dim\, \Sym^k(gl(m|n))$ for any $k$ and generic $q$ (see \cite{GPS2}). Thus, we can  define the product in the
algebra $\lq$ as a new noncommutative product in the initial (super-commutative) algebra $\Sym(gl(m|n))$ so that this product depends smoothly on $q$.

The explicit construction can be shortly described as follows. In \cite{GPS2} we discussed the projectors $P^k:\L^{\ot k}\to \L^{\ot k}$
(called {\em $R$-symmetrizers}). Explicitly they are known only for $k=2,3$. Define the map $\al_k:\Sym^k(gl(m|n))\to \lq^k$ as follows.
Let us embed  $\Sym^k(gl(m|n))$ into $\L^{\ot k}$ in a natural way (we identify
$\L$ and $gl(m|n)$ as linear spaces). Then $\al_k$ is by definition the restriction of $P^k$ to $\Sym^k(gl(m|n))$. Using the
family of maps $\{\al_k\}$ we can push the product in the algebra $\lq$ to that in $\Sym(gl(m|n))$ by
$$
f\star_q g=\al^{-1}_{k+l}(\al_k(f)\al_l(g)),\quad {\rm if} \quad  f\in \Sym^k(gl(m|n)),\; g\in \Sym^l(gl(m|n)).
$$

Now develop this product in $\nu=\log(q)$
$$
f\star_q g=fg+\nu c_1(f,g)+\nu^2 c_2(f,g)+..., c_i(f,g)\in \Sym(gl(m|n)).
$$
Following  the classical pattern (with one additional condition indicated below) we can show that the expression
$\{f,g\}=c_1(f,g)-c_1(\smn(f\ot g))$  is a super-Poisson bracket. This means that the following axioms are
fulfilled for $ \forall\,f,g,h \in \Sym(gl(m|n))$
$$
\begin{array}{c}
\{\,,\,\}=-\{\,,\,\} \smn (f\ot g),\\
\rule{0pt}{6mm}
\{\,,\,\}(\{\,,\,\}\ot I)(I+\sigma^{m,n}_{12}\sigma^{m,n}_{23}+ \sigma^{m,n}_{23}\sigma^{m,n}_{12})(f\ot g \ot h)=0, \\
\rule{0pt}{6mm}
\smn (\{f,g\}\ot  h)=(I\ot\{\,,\,\})\sigma^{m,n}_{12}\sigma^{m,n}_{23}(f\ot g \ot h).
\end{array}
$$
Note that for  even or odd elements $f$ and $g$ the third relation leads to the consequence
$\overline{\rule{0pt}{3.5mm}\{f,g\}}=\bar{f}+\bar{g}$ where $\bar{f}$ is the parity of $f\in\Sym(gl(m|n))$.

The additional condition mentioned above consists in the following. The terms $c_1$ and $c_2$ should be coordinated with
the parity in the same manner:
$$
\smn(c_i(f,g)\ot h)=(I\ot c_i)\sigma^{m,n}_{12}\sigma^{m,n}_{23}(f\ot g \ot h),\quad \forall\,\, f,g,h \in \Sym(gl(m|n)),\, i=1,2.
$$

Assuming the Hecke symmetry to be of the form (\ref{R-matr}), we conjecture that  the product $\star_q$  in the algebra $\Sym(gl(m|n))$ is such that the corresponding terms $c_1$ and $c_2$ possess this property. (In order to check this
conjecture we need an explicit form of the $R$-symmetrizers $P^k$ mentioned above.)

Taking this conjecture for granted, it is not difficult to compute the corresponding Poisson bracket. Fist, we rewrite the
relations (\ref{RE}) (with $\hbar= 0$) as follows
$$
\RR L_1 \RR_{21} L_2-L_2 \RR L_1 \RR_{21} = 0,
$$
where $\RR=\smn R$, $\RR_{21}= R\, \smn$ and $L_2=\smn L_1 \smn$.
Then by developing  the operator $\RR=I+\nu r+...$ and comparing the terms linear in $\nu$, we find the corresponding
Poisson bracket on the generators of the algebra $\Sym(gl(m|n))$
\be
\{\,,\,\}'(L_1\ot L_2)=-r L_1 L_2-L_1 r_{21} L_2+L_2 L_1 r_{21}+L_2 r L_1. \label{bra}
\ee
Here $r_{21}=\smn r \smn$ and $r\in \End(\vv)$ is given by the following formula
$$r_{ij}^{kl}=\sum_{1\leq i\leq m+n} (-1)^{\oi} \,(1-2\oi) e_{i}^i\ot e_{j}^j+2\sum_{j>i} e_{j}^i\ot e_{i}^j. $$

Similarly to the classical case this Poisson bracket is compatible with the linear bracket $\{\,,\}_{gl(m|n)}$.
These two bracket span the pencil which is the semi-classical counterpart of  the algebra $\lqh$.

Besides,  all elements $\Tr_\smn L^k$ belong to the Poisson center of any bracket from the Poisson pencil spanned by the
brackets $\{\,,\,\}_{gl(m|n)}$ and $\{\,,\,\}'$.
In other words, this pencil has the center of $gl$ type.
So, any such a bracket restricts to all "super-orbits" defined by polynomial
equations $\Tr_\smn L^k=a_k$. These Poisson structures and their braided analogs are considered in the next section.

\section{Poisson pencils on super-orbits and their quantization}
\label{sec:4}

As was shown in \cite{GPS1}, if $R$ is any skew-invertible Hecke symmetry, then the generating matrix $L$ of the algebra
$\lq$ satisfies a Cayley-Hamilton (CH) type identity. If $R$ is a deformation of the super-flip $\smn$ this identity takes the
form
\be
\sum_{i=1}^{m+n} b_i(L) L^i=0
\label{CH}
\ee
where $b_i(L)$ are non-trivial central elements of the algebra $\lq$. Moreover, if $n\not=0$ the leading coefficient
$b_{m+n}(L)$ is not a number.

Upon multiplying the identity (\ref{CH}) by $b_{m+n}(L)$, we can represent it in the following factorized form
$$
\big(b_{m+n}(L)\prod_{i=1}^m (L-\mu_i)\big) \,\big(b_{m+n}(L)\prod_{j=1}^n (L-\nu_j)\big)=0.
$$
Here the {\em eigenvalues} $\mu_i$ and $\nu_j$ are elements of the algebraical extension of $Z(\lq)_{loc}$ where
$Z(\lq)_{loc}$ stands for the localization of the center $Z(\lq)$  by the set
$\{b_{m+n}^k(L),\, k\ge 1\}$. The eigenvalues $\mu_i$ (resp., $\nu_j$) are called {\em even} (resp., {\em odd}).

The reason for this terminology is the following formula expressing the quantities $\Tr_\smn L^{k}$ via these eigenvalues
for a super-matrix:
$$
\Tr_\smn L^{k}=\sum_{i=1}^m \mu_i^{k}-\sum_{j=1}^n \nu_j^{k},\quad k\ge 0.
$$
Below we give a braided analog of this formula (see (\ref{pk}), (\ref{dd-i})).

In what follows we consider the quotient algebras
\be
\K[\O_{\mu,\nu}]=\Sym(gl(m|n))/ J_{\smn}(\mu,\nu)
\label{sup-orb}
\ee
where the ideal $J_{\smn}(\mu,\nu)$ is generated by $m+n$ elements
$$
\Tr_{\smn}L^k - \Big(\sum_{i=1}^m\mu_i^k - \sum_{j=1}^n\nu_j^k\Big),\quad
1\le k\le m+n.
$$
The algebras $\K[\O_{\mu,\nu}]$ are super-analogs of (the coordinate algebras of) affine algebraical varieties.

Let  $\{\,,\,\}_{a,b}^{\O_{\mu,\nu}}$ be the restriction  of the bracket $\{\,,\,\}_{a,b}$ to this super-variety. Now, we want
to discuss two questions.

For what values of $\mu=(\mu_1,\mu_2,...,\mu_m)$ and $\nu=(\nu_1,\nu_2,...,\nu_n)$ the quotient $\K[\O_{\mu,\nu}]$ can
be considered as a {\em regular} super-variety and how to
quantize the pencil $\{\,,\,\}_{a,b}^{\O_{\mu,\nu}}$? Hopefully, a quantum analog of $\K[\O_{\mu,\nu}]$ is a {\em braided
variety} which can be presented in a similar way
\be
\lq/\langle \Tr_R L-a_1,\, \Tr_R L^2-a_2, ..., \Tr_R L^{m+n}-a_{m+n}  \rangle,\quad a_i\in \K.
\label{vari}
\ee
However, first we should answer an analogous question in the quantum case: for which values of numbers $a_i$
the quotient (\ref{vari}) of the quantum algebra $\lq$ can be considered as a {\em regular} braided variety?
We have to answer this question since it is natural to suppose that a quantum counterpart of a  regular (super-)variety is a regular braided  one.

It is known that in a classical case $(n=0, q=1)$ the variety $\K[\O_\mu]$ is regular iff it is a generic orbit, i.e. the orbit of
a matrix with pairwise distinct eigenvalues $\mu_i$. If an affine algebraical  variety $\M$ is defined by a system of polynomial
equations $p_i=0,\,\, 1 \leq i\leq k$ then it is regular iff the rank of a matrix formed by gradients of $p_i$ is maximal at each
point of the variety. If it is the case, then according to the known Serre result (see \cite{S}) the space of sections of the
cotangent bundle on $\M$ is a finitely generated projective $\K[\M]$-module. We call this module  {\em cotangent}.

In \cite{GS2} we succeeded in constructing analogs of the cotangent module over super- and braided varieties
((\ref{sup-orb}) and (\ref{vari}) respectively) for generic values of the quantities $\Tr_R L^k$. However, construction of such a
module fails for some {\em exceptional} values of these quantities. We want  to describe the set of exceptional values in terms
of eigenvalues of the matrix $L$. To this end we employ the formula expressing $\Tr_R L^k$ in terms of $\mu_i$ and $\nu_j$
(see \cite{GPS3}):
\begin{equation}
\Tr_R L^k=\sum_{i=1}^m d_i \mu_i^k +
\sum_{j=1}^n  d'_j\nu_j^k\quad \forall\,k\ge 0\,,
\label{pk}
\end{equation}
where the {\em quantum dimensions} $d_i$ and $d'_j$ read
\be
d_i = q^{-1}\prod_{p=1\atop p\not=i}^m
\frac{\mu_i - q^{-2}\mu_p}{\mu_i-\mu_p}\,
\prod_{j=1}^n \frac{\mu_i - q^2\nu_j}{\mu_i -\nu_j} , \quad
d'_j = -\,q\,\prod_{i=1}^m \frac{\nu_j - q^{-2}\mu_i}{\nu_j-\mu_i}\,
\prod_{p=1\atop p\not=j}^n \frac{\nu_j - q^2\nu_p}{\nu_j-\nu_p}\,.
\label{dd-i}
\ee

Thus, expressing the coefficients $a_i$ in (\ref{vari}) in terms of the eigenvalues $(\mu, \nu)\in \K^{\oplus(m+n)}$ we present
the algebra (\ref{vari}) as the quotient
\be
\Omn=\lq/ J_{R}(\mu,\nu),
\label{vari1}
\ee
where the ideal $J_R(\mu,\nu)$ is generated by the following elements
\be
\Tr_R L^k-\Big(\sum_{i=1}^m d_i \mu_i^k + \sum_{j=1}^n d'_j\nu_j^k\Big),
\quad 1\le k\le m+n.  \label{JR}
\ee
with $d_i$ and $d'_j$ given by (\ref{dd-i}).

As was shown in \cite{GS2}, the cotangent module exists on such a braided variety iff the following conditions are fulfilled
\be
\mu_i\not=q^2\mu_j,\quad \nu_i\not=q^2\nu_j,\quad \mu_i\not=q^2\nu_j,\, \quad  1\leq i \leq m, 1\leq j \leq n.
\label{cond}
\ee
Let $\E$ be the set of eigenvalues $(\mu, \nu)\in \K^{\oplus(m+n)}$ such that at least one of these conditions fails.
We call the algebra $\Omn$ with $(\mu, \nu)\in \K^{\oplus(m+n)}\backslash \E$ a {\em braided generic orbit}.

Under the limit $q\to 1$ we get a similar condition for a {\em
generic super-orbit}. In this case an analog of the restrictions
(\ref{cond}) reads \be \mu_i\not=\mu_j,\quad \nu_i\not=\nu_j,\quad
\mu_i\not=\nu_j,\, \quad  1\leq \leq m, 1\leq j \leq n.
\label{cond1} \ee Thus, as a quantization of the Poisson bracket
$\{\,,\,\}'$ restricted to a generic super-orbit $\K[\O_{\mu,
\nu}]$ we can consider the braided variety $\Omn$ with the same
eigenvalues $(\mu, \nu)$. It is evident that if $q-1$ is small
enough, then the conditions (\ref{cond1}) entails these
(\ref{cond}). Consequently, the corresponding braided variety
$\K_q[\O_{\mu, \nu}]$ is regular or, in other words, a braided
generic orbit.  However, another choice of the quantum object is
also possible: we have only to verify the conditions (\ref{cond}).

In conclusion of this section we want to emphasize that {\em non-commuta\-tive}  super-(or brai\-ded) varieties can be
considered in a similar manner. They are appropriate quotients of the modified REA $\lqh$. Namely, they are defined by the
same formula (\ref{vari1}) but with quantum dimensions given by
\be
\begin{array}{l}
\displaystyle
d_i = q\,\prod_{p=1\atop p\not=i}^m
\frac{\mu_i - \qqq\mu_p-\qq\h}{\mu_i-\mu_p}\,
\prod_{j=1}^n \frac{\mu_i - q^2\nu_j+q\h}{\mu_i -\nu_j}\, ,\\
\rule{0pt}{10mm}
\displaystyle
d'_j = -\, q \,\prod_{i=1}^m \frac{\nu_j - \qqq\mu_i-\qq\h}{\nu_j-\mu_i}\,
\prod_{p=1\atop p\not=j}^n \frac{\nu_j - q^2\nu_p+q\h}{\nu_j-\nu_p}\,
\end{array}
\label{d-i}
\ee
(see \cite{GS2}).

Consequently, the conditions (\ref{cond}) must be modified as well. Thus, the quotient of the algebra $\lqh$ is by definition
a {\em regular braided non-commutative orbit} iff
$$
\mu_i - q^{-2}\mu_j-\qq\h\not=0,\, \nu_j - q^2\nu_j+q\h\not=0,\,\mu_i - q^2\nu_j+q\h\not=0.
$$
Similarly to regular braided varieties considered above, this definition is motivated by the fact that on such a regular
braided non-commutative orbit there exists the cotangent module. For detail the reader is refereed to \cite{GS2}.

Finally, the quotient $\lqh/\hat{J}_R(\mu,\nu)$, where the ideal $\hat{J}_R(\mu,\nu)$ is defined by the formula similar to that (\ref{JR}) but
with $d_i$ and $d'_j$ given by formulae (\ref{d-i}), is just a quantum counterpart of the pencil spanned by the brackets $\{\,,\,\}_{gl(m|n)}$ and
(\ref{bra}) restricted to the super-orbit $\Omn$.

\section{Examples}
\label{sec:5}

Let us complete consideration of the example of section \ref{sec:2}. Here we treat it from the viewpoint developed in
sections \ref{sec:3} and \ref{sec:4}.

Consider the Hecke symmetry (\ref{R-matr}) for $m=2$, $n=0$
$$
R = \left(
\begin{array}{cccc}
q&0&0&0\\
0&q-q^{-1}&1&0\\
0&1&0&0\\
0&0&0&q\end{array}\right).
$$

The defining relations of the corresponding algebra $\lqh$ with the generating matrix
$L=\left(
\begin{array}{cc} a&b\\
c&d\end{array}\right)$ are as follows
$$
\begin{array}{c}
q ab-\qq ba=\h b,\quad q ca-\qq ac=\h c,\quad ad-da=0,\\
\rule{0pt}{6mm}
q(bc-cb)=(\la a-\h)(d-a),\quad q(cd-dc)=c(\la a -\h),\quad q(db-bd)=(\la a -\h)b.
\end{array}
$$

The operators $B$ and $C$ are given by the matrices
$$
B= \left(
\begin{array}{cc}
\qq&0\\
0&q^{-3}\end{array}\right),\quad
C= \left(
\begin{array}{cc}
q^{-3}&0\\
0&\qq\end{array}\right).
$$

Thus, we have
$$
\ell=\Tr_R L=q^{-3} a+ \qq d,\quad \Tr_R L^2= q^{-3} (a^2+bc)+\qq (cb+d^2).$$
These elements are central in the algebra $\lqh$. For the numerical $R$-trace $\tr_R: \L\to \K$ on the space $\L={\span}_\K(a,b,c,d,)\cong \End(V)$
we get
$$
\tr_R\, a=1,\quad \tr_R\, b=0,\quad \tr_R\, c=0,\quad \tr_R\, d=1.
$$
Therefore, the elements $b, c, h=a-d$ are traceless. Besides, the pairing (\ref{pai}) takes the form
$$
\langle a,a \rangle=q,\quad \langle d,d \rangle=\qq,\quad \langle b,c \rangle=\qq,\quad \langle c,b \rangle=q,
$$
all other terms being zero. The first and second formulae above are equivalent to
$$
\langle h,h \rangle=2_q, \quad
\langle \ell,\ell \rangle=q^{-4} 2_q.
$$

On rewriting the defining relations for $\lqh$ in the basis $\ell, b,c, h$ and setting $\ell = 0$,
we recover the defining relations of the algebra $\slqh=\lqh/\langle \ell \rangle$:
\be
q^2 hb-bh=2_q\h b,\quad q^2 ch-ch= 2_q\h c,\quad (q^2+1)(bc-cb)+(q^2-1) h^2= 2_q\h h.
\label{motiv}
\ee
They coincide (up to a notation) with relations (\ref{quant2}) whith $A=B=C=2_q\h$.
The quadratic central element $\Tr_R L^2$ being reduced to the algebra $\slqh$ becomes
$\Cas_q$ as in section \ref{sec:2} where it  was found by other means. Being restricted to the
space $\SL={\span}_\K(b,h,c)$, the pairing (\ref{pai}) takes the form (\ref{pair}) (up to a factor).

The CH identity for the matrix $L$ reads:
$$
L^2-(q^{-2}a+d) L+(q^{-2}ad - cb)I=0.
$$
Thus, according to our definition of eigenvalues we have
$$
\mu_1+\mu_2=q^{-2}a+d,\quad \mu_1\mu_2= q^{-2}ad - cb.
$$
Expressing the quantities $\Tr_R L$ and $\Tr_R L^2$ via these eigenvalues we introduce
a {\em braided variety} by the following system of polynomial equations
$$
\begin{array}{c}
 \Tr_R L=q^{-3}a+\qq d=\qq (\mu_1+\mu_2),\\
\rule{0pt}{6mm}
\Tr_R L^2=q^{-3}(a^2+bc)+ \qq(cb+d^2)=\qq (\mu_1^2+\mu_2^2)+(\qq-q^{-3}) \mu_1\mu_2.
\end{array}
$$
Such a variety is a braided {\em generic} orbit iff $\mu_1\not=q^{\pm 2} \mu_2$. By imposing the
condition $Tr_R L=0$ we get the braided analog of a hyperboloid. This condition entails $\mu_2=-\mu_1$.
So, the {\em braided hyperboloid} can be parameterized by one parameter, for instance by $\mu_1$. Explicitly, it is given by the following equation
$$
\Tr_R L^2=(q^{-3}+\qq) \mu_1^2.
$$
For a generic $q$ it is a regular braided variety (and consequently generic orbit) for any $\mu_1\not=0$.

Braided non-commutative orbits can be defined in a similar way as appropriate quotients of the algebra $\lqh$. For this
purpose we have to replace the above system of equations by the following one
$$
\begin{array}{c}
\Tr_R L=q^{-3}a+\qq d=\qq (\mu_1+\mu_2)-\qqq\h,\\
\rule{0pt}{6mm}
\Tr_R L^2=q^{-3}(a^2+bc)+ \qq(cb+d^2)=\qq (\mu_1^2+\mu_2^2)+(\qq-q^{-3}) \mu_1\mu_2-\qqq\h (\mu_1+\mu_2).
\end{array}
$$
This braided  non-commutative variety is a non-commutative generic orbit iff
$$\mu_1\not=\qqq \mu_2-\qq\h \,\,{\rm  and}\,\, \mu_2\not=\qqq \mu_1-\qq\h.$$

To obtain the corresponding bracket $\{\,,\,\}'$ on the space $gl(2)$ it suffices to extend the bracket constructed in
section \ref{sec:2} by the generator $\ell$ which is Poisson commuting  with other generators.

Now, consider another example related to the super-Lie algebra  $gl(1|1)$. This algebra is generated by 4 elements
$a,b,c,d$ subject to the  relations
$$
\begin{array}{c}
[a,b]=b,\quad [a,c]=-c,\quad [a,d]=0, \quad [d,b]=b,\\
\rule{0pt}{6mm}
[d,c]=-c,\quad [b,c]_+=d-a,\quad [b,b]_+=[c,c]_+=0.
\end{array}
$$
Notation $[\,,\,]_+$ stands for the anti-commutator. Emphasize that this basis differs from the usual one by the sign
at $b$ and $d$, our choice is motivated by that in the algebra $\lqh$. The elements
\be
\Tr_{\sigma^{1,1}} L=a-d\quad {\rm  and}\quad \Tr_{\sigma^{1,1}} L^2=a^2+bc-cb-d^2
\label{sledy}
\ee
are central in the enveloping algebra $U(gl(1|1))$.

The corresponding Hecke symmetry is
$$
R = \left(
\begin{array}{cccc}
q&0&0&0\\
0&q-q^{-1}&1&0\\
0&1&0&0\\
0&0&0&-q^{-1}\end{array}\right).
$$

The operators $B$ and $C$ are
$$
B= \left(
\begin{array}{cc}
\qq&0\\
0&-q^{-1}\end{array}\right),\qquad
C= \left(
\begin{array}{cc}
q&0\\
0&-q\end{array}\right).
$$

The defining relations in corresponding  algebra $\lqh$  are as follows
$$
\begin{array}{l}
q^2 ab-ba=q\h b,\quad q^2 ca-ac=q\h c,\quad ad-da=0,\quad b^2=c^2=0,\\
\rule{0pt}{6mm}
\qq bc+qcb -(q-\qq) a(a-d)=\h (a-d),\\
\rule{0pt}{6mm}
bd-db-(q^2-1) ab=q \h b,\\
\rule{0pt}{6mm}
cd-dc +(q^2-1) ca= -q\h c.
\end{array}
$$

The related bracket $\{\,,\,\}'$ is
$$
\begin{array}{l}
\{a,b\}'=ab,\quad \{a,c\}'=-ac,\quad \{a,d\}'=0,\\
\rule{0pt}{6mm}
\{d,b\}'=ab,\quad {\{b,b\}'}_+={\{c,c\}'}_+=0,\\
\rule{0pt}{6mm}
\{d,c\}'=-ac,\quad {\{b,c\}'}_+=cb-a(a-d).
\end{array}
$$
Note that the elements $b,c,h=a-d$ generate the super-Lie subalgebra $sl(1|1)$.
Its multiplication table is
\be
[h,b]=[h,c]=0,\,\, [b,c]_+=-h.
\label{rr}
\ee
However,  the  bracket $\{\,,\,\}'$ has no restriction to the super-algebra $\Sym(sl(1|1))$. Also, note that the elements
(\ref{sledy}) are central for any bracket from the pencil spanned by the brackets $\{\,,\,\}_{gl(1|1)}$ and $\{\,,\,\}'$.

In this case the CH identity for the matrix $L$ takes the from
$$
(a-d)L^2-(a^2-d^2+bc-cb)L+((a-d)(bc-ad)-a(bc-cb))I=0.
$$
This identity can be written in the factorized form (after additional multiplication by $\ell = \Tr_RL = q(a-d)$)
$$
\Big(\ell L-q S(L)\Big)\Big(\ell L+q^{-1}A(L)\Big) = 0,
$$
where the polynomials $S(L)$ and $A(L)$ read
\be
S(L) = \frac{1}{2_q}\,\Big(q^{-1}\ell^2 + \Tr_RL^2\Big),\qquad
A(L) = \frac{1}{2_q}\,\Big(q\,\ell^2  - \Tr_RL^2\Big).
\label{sh-f}
\ee
Thus, the even $\mu$ and odd $\nu$ eigenvalues of $L$ are defined by the fractions
\be
\mu = q\,\frac{S(L)}{\ell},\qquad
\nu = -q^{-1}\frac{A(L)}{\ell}.
\label{eig}
\ee

Relations (\ref{sh-f}) and (\ref{eig}) allow us to express $\Tr_RL$ and $\Tr_RL^2$ in terms of
eigenvalues. Thus, according to our general approach a {\em braided variety} is defined in the algebra $\lq$
by the system of equations
$$
\begin{array}{c}
\Tr_R L=q (a-d)=\qq \mu- q\nu,\\
\rule{0pt}{6mm}
\Tr_R L^2=q(a^2+bc-cb-d^2)=(\mu+\nu)(\qq \mu- q\nu).
\end{array}
$$
Emphasize that the braided variety defined by the equation
$$
\Tr_R L=q (a-d)=\qq \mu- q\nu=0
$$
is not regular and according to our terminology is not a braided generic  orbit.

Turning to non-commutative braided varieties we have to replace the above equations by
$$
\begin{array}{c}
\Tr_R L=q (a-d)=\qq \mu- q\nu+\h,\\
\rule{0pt}{6mm}
\Tr_R L^2=q(a^2+bc-cb-d^2)=(\mu+\nu)(\qq \mu- q\nu)+\h(\mu+\nu).
\end{array}
$$
The non-commutative braided varieties defined in the algebra $\lqh$ by these equations are regular
iff $\mu\not=q^2\nu-q\h$.

\section{Open problems and concluding remarks}
\label{sec:6}

Here we formulate some open problems. Consider a generalization of
the bracket $\{\,,\,\}'$ (see section \ref{sec:2}) realized in
terms of the compact form of $sl(2,\C)$. Namely, let $\{x,y, z\}$ be
the standard basis in the polynomial algebra $\K[so(3)^*]\cong
\Sym(so(3))$. In this algebra we introduce a Poisson bracket defined on the generators as follows
 \be \{x,y\}'=z\, p(x, y,z)
,\quad \{y,z\}'=x\, p(x, y,z),\quad \{z,x\}'=y\, p(x, y,z),
\label{comp} \ee where $p(x, y,z)$ is a fixed polynomial in
$x,y,z$. We leave checking the fact that all axioms of a Poisson bracket are fulfilled to the reader.

Also, observe that this Poisson bracket is compatible with that
$\{\,,\,\}_{so(3)}$. Observe that
the center of each bracket from the  pencil spanned by these
$\{\,,\,\}_{so(3)}$ and (\ref{comp})
consists of the elements
$q(x^2+y^2+z^2)$ where $q$ is a polynomial in one variable.

 The following questions are of great interest.

\begin{enumerate}
\item
How to  explicitly quantize these brackets?

\item
What is the center of a quantum algebra obtained by a quantization, what is the corresponding
pairing, and whether it is possible to treat this pairing via a deformed trace?

\item
How to quantize these brackets restricted to a sphere and whether it is possible to describe the
corresponding quantum algebras in the same way as braided varieties above?

\item
How to classify all Poisson structures on $gl(m)^*$ (or $gl(m|n)^*$) possessing the same center as
the Poisson-Lie bracket $\{\,,\,\}_{gl(m)}$ (or $\{\,,\,\}_{gl(m|n)}$) has and how to quantize them as well
as their restrictions to a generic orbit in $gl(m)^*$ (or $gl(m|n)^*$)?
\end{enumerate}

Anyway, the corresponding quantum objects can not be described in the
frameworks of the above braided geometry in its present limits which must be extended.

{\bf Remark 1.} In conclusion we would like to stress that the above Poisson structures on the
space $gl(m)^*$ (apart from the Poisson-Lie one) are not unimodular. Roughly speaking, a Poisson
structure is called unimodular if there exists a volume form in a sense compatible with the defining  bi-vector field.
A particular case of unimodular Poisson structures are simplectic ones defined via closed 2-forms. For a simplectic
Poisson structure the role of such a  volume form is played by the Liouville measure $\Om$. For this measure the following relation is valid
$$
\int \{f,g\}\Om=0,\,\,\,\, \forall f,\, g .
$$
So,  the map $f\to \int f \Om$ is a Poisson analog of the trace. A quantization of such a Poisson structure gives rise to
an algebra with a trace possessing the classical  property $\Tr[.,.]=0$ (see \cite{GR}).

A standard  example of a simplectic Poisson structure is the restriction of the bracket $\{\,,\,\}_{so(3)}$ to a sphere
$\Cas=r^2$. Its quantum analog is a proper quotient of the algebra $U(so(3))$  (we put $\h=1$). For some discrete values
of $r$ the quantum algebra can be represented in a finite dimensional Hilbert space endowed with the usual trace. However,
the restriction to this sphere of the bracket (\ref{comp}) with $p(x,y,z)=z$ is not simplectic. Its Poisson leaves are two
half-spheres $z>0$, $z<0$ and each point of the equator $z=0$. In general, the brackets from the corresponding Poisson
pencil $\{\,,\,\}_{a,b}$ are not simplectic either. Nevertheless, their quantum counterparts which are appropriate quotients
of the algebras $\L(1,q)$ can be represented in  finite dimensional spaces (also for some special values of eigenvalues $\mu,
\nu$). The essential point here is that the usual trace should be changed for a deformed (braided) trace. A category of
such representation is considered in \cite{GPS2}. Thus, though a Poisson analog of the trace does not exist on the whole
sphere the corresponding quantum algebra can be endowed with a trace but this trace is {\em braided}.
\medskip

{\bf Remark 2.} We would like to emphasize a difference between a usual variety and a super-one. Consider again
the algebra $sl(1|1)$. The quotient
$$
\Sym(sl(1|1))/\langle h^2+bc-r^2\rangle,\quad r\in\K
$$
is a regular super-variety  iff $r\not=0$. A construction of the cotangent module (which is projective) in the case  $r\not=0$ is evident.

However, in virtue of the factorization
$$
h^2+bc-r^2=(h-(r-{{bc}\over{2r}}))(h+(r-{{bc}\over{2r}}))
$$
the super-variety in question is a union of two super-varieties: each of them is defined by one of the equations
\be
h-(r-{{bc}\over{2r}})=0,\qquad h+(r-{{bc}\over{2r}})=0.
\label{ss}
\ee

Nevertheless, these two super-varieties have no common points. Indeed, the system (\ref{ss})
is equivalent to that $h=0$, $bc=2r^2$ where the second equation has no solution if $r\not=0$.
By contrast, the  system (\ref{ss}) in the classical case (i.e. if  all generators  are even)
describes a non-empty  set of points. All these points of the variety defined by the equation
\be
(h-(r-{{bc}\over{2r}}))(h+(r-{{bc}\over{2r}}))=0 \label{last}
\ee
are singular. Thus, the  variety (\ref{last}) is regular or not in function of the parity of the generators $b$ and $c$.
Hopefully, braided deformations of super-orbits above can be presented in a form similar to that (\ref{last}).

\end{document}